\documentclass[10pt]{article} \usepackage{latexsym}  
\usepackage{amsfonts} 
\usepackage{amsmath} 
\usepackage{amssymb} 
\usepackage{amscd}   
\usepackage{multicol}    
\usepackage{enumerate} 
\usepackage{amsthm}            
\theoremstyle{plain}
    \newtheorem{theorem}                    {Theorem}       [section]
    \newtheorem{lemma}      [theorem]       {Lemma}
    \newtheorem{corollary}  [theorem]       {Corollary}
    \newtheorem{proposition}[theorem]       {Proposition}

\begin{document}

\newcommand{\cosk}{\operatorname{cosk}}
\newcommand{\sk}{\operatorname{sk}}
\newcommand{\C}{\mathcal C}
\newcommand{\zz}{\mathfrak z}
\newcommand{\sh}{\mathcal S}
\newcommand{\chr}{\operatorname{char}}
\newcommand{\Aut}{\operatorname{Aut}}
\newcommand{\Hom}{\operatorname{Hom}}
\newcommand{\Tor}{\operatorname{Tor}}
\newcommand{\End}{\operatorname{End}}
\newcommand{\Ext}{\operatorname{Ext}}
\newcommand{\Gal}{\operatorname{Gal}}
\newcommand{\Pic}{\operatorname{Pic}}
\newcommand{\ord}{\operatorname{ord}}
\newcommand{\tor}{{\operatorname{tor}}}
\newcommand{\Spec}{\operatorname{Spec}}
\newcommand{\op}{\mathrm{{op}}}
\newcommand{\holim}{\operatornamewithlimits{holim}}
\newcommand{\im}{\operatorname{im}}
\newcommand{\coim}{\operatorname{coim}}
\newcommand{\coker}{\operatorname{coker}}
\newcommand{\gr}{\operatorname{gr}}
\newcommand{\id}{\operatorname{id}}
\newcommand{\Br}{\operatorname{Br}}
\newcommand{\cd}{\operatorname{cd}}
\newcommand{\CH}{\operatorname{CH}}
\renewcommand{\lim}{\operatornamewithlimits{lim}}
\newcommand{\colim}{\operatornamewithlimits{colim}}
\newcommand{\rk}{\operatorname{rank}}
\newcommand{\codim}{\operatorname{codim}}
\newcommand{\NS}{\operatorname{NS}}
\newcommand{\N}{{\mathbb N}}
\newcommand{\Z}{{{\mathbb Z}}}
\newcommand{\Q}{{{\mathbb Q}}}
\newcommand{\R}{{{\mathbb R}}}
\newcommand{\F}{{{\mathbb F}}}
\newcommand{\f}{{\mathcal F}}
\renewcommand{\P}{{{\mathcal P}}}
\newcommand{\OO}{{{\mathcal O}}}
\renewcommand{\L}{{\mathcal L}}
\newcommand{\A}{{\mathcal A}}
\newcommand{\Sm}{\text{\rm Sm}}
\newcommand{\Sch}{\text{\rm Sch}}
\newcommand{\et}{{\text{\rm\'et}}}
\newcommand{\Zar}{{\text{\rm Zar}}}
\newcommand{\Nis}{{\text{\rm Nis}}}
\newcommand{\cdh}{\text{{\rm cdh}}}
\newcommand{\Div}{\operatorname{Div}}
\newcommand{\Alb}{\operatorname{Alb}}
\renewcommand{\div}{\operatorname{div}}
\newcommand{\corank}{\operatorname{corank}}
\renewcommand{\O}{{\mathcal O}}
\newcommand{\p}{{\mathfrak p}}
\newcommand{\pschemes}{{\operatorname{Sch}^\bullet_k}}
\newcommand{\proofend}{\hfill $\square$ \\ \smallskip}

\title{Homological descent for motivic homology theories}
\author{Thomas \textsc{Geisser}\footnote{Graduate School of Mathematics,
Nagoya University, Furucho, Nagoya 
464-8602, Japan.\newline e-mail: \texttt{geisser@math.nagoya-u.ac.jp}}}

\date{}
\maketitle

\begin{abstract}      
The purpose of this paper is to give homological descent theorems for 
motivic homology theories (for example Suslin homology) 
and motivic Borel-Moore homology theories (for example higher Chow groups)
for certain hypercoverings.
\end{abstract}

\section{Introduction} 
We consider covariant 
functors $\f$ from the category of schemes separated and of 
finite type over a fixed noetherian scheme $X$ and proper morphisms 
to the category of homologically positive complexes of 
abelian groups. 
We assume that 
for every abstract blow-up square, there is a distinguished triangle
\begin{equation}\label{blowup}
\f (Z')\to \f (X')\oplus \f (Z) \to \f (X) \to \f (Z')[1]
\end{equation}
in the derived category of abelian groups.
Here an abstract blow-up square is a diagram of the form 
\begin{equation}\label{blbl} 
\begin{CD}
Z'@>>> X'\\
@VVV @V\pi VV\\
Z@>i>> X
\end{CD}
\end{equation}
with $\pi$ proper, $i$ a closed embedding, and such that $\pi$
induces an isomorphism $X'-Z' \to X-Z$. The morphism 
$Z\coprod X'\to X$
is then called an abstract blow-up. In particular, (taking $X'=\emptyset$),
every closed embedding $Z\to X$ defined 
by a nil-potent ideal indudes a quasi-isomoprhism $\f(X')\to \f(X)$. 

Recall that a proper cdh-cover is a proper map $p:X'\to X$ such for 
every point
$x\in X$, there is a point $x'\in X'$ with $p(x')=x$ and 
$p^*:k(x)\stackrel{\sim}{\to} k(x')$,
and a hyperenvelope is an augmented simplicial scheme $a:X_\bullet \to X$
such that for every $n$, the map $X_{n+1} \to (\cosk_n X_\bullet)_{n+1}$
is a proper cdh-cover.
For a simplicial scheme $X_\bullet$, we simply write $\f(X_\bullet)$
for the total complex of the simplicial complex of abelian groups.

\begin{theorem}\label{main1}
For any functor as above,  
and for any hyperenvelope $a:X_\bullet\to X$, 
the augmentation map induces a quasi-isomorphism
$$ \f (X_\bullet)\to \f (X).$$
\end{theorem}

In characteristic $p$, smooth hyperenvelopes are only known to exist under
resolution of singularities. To remove this hypothesis, and to be able to 
use Gabber's refinement of de Jong's theorem on alterations,
we consider $l$-hyperenvelopes. A proper
ldh-cover is a proper surjection $p:X'\to X$ such that for every point 
$x\in X$ there is a point $x'\in X'$ with $p(x')=x$ and such 
$k(x')$ is a finite extension of degree prime to $l$ of $p^*(k(x))$.
An $l$-hyperenvelope is an augmented simplicial scheme $a:X_\bullet\to X$
such that $X_{n+1} \to (\cosk_n X_\bullet)_{n+1}$ is a proper ldh-cover
for all $n$.


\begin{theorem}\label{main2}
Assume that $\f $ is a functor to the category of complexes of 
$\Z_{(l)}$-modules, which satisfies in addition to the above
the following property:

For any finite flat map $p:X\to Y$ of degree $d$ prime to $l$, 
there is a functorial pull-back map $p^*:\f (Y)\to \f (X)$ 
such that $p_*p^*$ induces multiplication by $d$ on homology, and which 
is compatible with base-change.

Then for any $l$-hyperenvelope $a:X_\bullet\to X$, 
the augmentation map induces a quasi-isomorphism  
$$ \f (X_\bullet)\to \f (X).$$
\end{theorem}

One can see that some hypothesis on the coefficients is necessary by considering the 
\v Cech-nerve $\cosk_0(L/k)$ 
of a finite field extension $L/k$ of degree $d$: 
The map $H_0^S(L,\Z)\to H_0^S(k,\Z)$ is multiplication by $d$ on $\Z$, 
hence descent does not hold. 

The proof of the theorems is along the lines SGA 4 Vbis \S  3. 
Gillet \cite{gillet} used a similar argument 
to prove descent for higher Chow groups and $K'$-homology, 
but, using the notes of B.Conrad \cite{brian},
we give a self-contained proof 
which in addition does not require the localization property.


As an application, we obtain the following descent theorem for the
motivic homology groups $H_{i}(X,A(n))$, motivic Borel-Moore homology
groups  $H_i^{BM}(X,A(n))$, and higher Chow groups: 

\begin{theorem}\label{desccor}
Let $X$ be of finite type over a perfect field, and 
$A$ be an abelian group.  Suppose that $a:X_\bullet\to X$
is a hyperenvelope and resolution of singularities holds, 
or that $a:X_\bullet\to X$
is a $l$-hyperenvelope and that $A$ is a $\Z_{(l)}$-module. Then we have 
spectral sequences 
\begin{align*}
E^1_{p,q} = H_q(X_p,A(r)) &\Rightarrow H_{p+q}(X,A(r));\\
E^1_{p,q} = H_q^{BM}(X_p,A(r)) &\Rightarrow H_{p+q}^{BM}(X,A(r)).\\
\end{align*}
\end{theorem}

For a scheme essentially of finite type over a Dedekind ring,
we define higher Chow groups as the Zariski hypercohomology
of Bloch's cycle complex (for schemes over fields, 
this agrees with the homology of the global sectons of the cycle complex).

\begin{theorem}\label{desccorchow}
Let $X$ be of finite type over a field or a Dedekind ring, and $A$ be an abelian group.  
Suppose that $a:X_\bullet\to X$ is a hyperenvelope, 
or that $ A$ is a $\Z_{(l)}$-module and $a:X_\bullet\to X$
an $l$-hyperenvelope. Then we have spectral sequences for any $n\in \Z$,
$$E^1_{p,q} = CH_n(X_p,q,A) \Rightarrow CH_n(X,p+q,A).$$
\end{theorem}

The analogous results holds for $K'$-theory can be proven by the 
same method.

We note that $l$-hyperenvelopes exist by by a theorem of Gabber:

\begin{theorem}\label{hyperexist}
For every scheme $U$ of finite type over a perfect field $k$ and 
any $l\not= \chr k$, 
there exists an $l$-hyperenvelope $U_\bullet\to U$ such that $U_\bullet$
is an open simplicial subscheme of a simplicial scheme consisting of 
smooth projective schemes over $k$.
\end{theorem}

\medskip
We thank S.Kelly for helpful comments, discussions and 
explanation of his work.

\section{Simplicial schemes and hyperenvelopes}
Let $\C$  be a category with finite limits. A simplicial object in 
$\C$ is a functor $X_\bullet: \Delta^\op\to \C$, where $\Delta$ is 
the simplicial category
of finite ordered sets $[i]=\{ 0,\ldots ,i\}$ with non-decreasing maps.
As usual, we write $X_n$ for $X_\bullet([n])$ and $\alpha^*:X_j\to X_i$
instead of $X_\bullet(\alpha)$ for $\alpha: [i]\to [j]$.
If $\Delta_{\leq n}$ is the full subcategory
of $\Delta$ consisting of $[0],\ldots , [n]$, then the restriction 
functor $i_n^*$ from simplicial objects to restricted simplicial objects,
i.e. functors $\Delta^\op_{\leq n}\to \C$, has a left adjoint skeleton 
$\sk_n$ and a right 
adjoint coskeleton $(i_n)_*= \cosk_n$. We note that in the literature, 
the notation
$ \sk_n$ both appears as the name of the restriction functor 
\cite{deligne} as well as the name of its left adjoint (e.g. SGA 4 V 7). 
By abuse of notation, we also denote the composition of
$(i_n)_*i_n^*$ by $\cosk_n$. In this notation,
the adjunction map takes the form $X_\bullet \to \cosk_nX_\bullet$. 
Concretely, 
$$(\cosk_nX_\bullet)_m = \lim_{D_m'}X_\phi$$
where $D_m'$ is the category of non-decreasing maps
$\phi:[i]\to [m]$ with $i\leq n$, $X_\phi= X_i$, 
and morphisms $\alpha: \phi\to \phi'$ the maps 
$[i]\to [i']$ compatible with the maps to $[m]$. 
This can also be expressed as follows \cite[Cor.3.10]{brian}:
Let $D_m$ be the full subcategory of $D_m'$ with objects increasing 
{\it injections} $\phi: [i]\to [m]$ for $i\leq n$ (which implies
that morphisms $\alpha: \phi\to \phi'$  are also is injective). Then  
$(\cosk_nX_\bullet)_m=X_m $ for $m\leq n$ and, for $m>n$,
$ (\cosk_nX_\bullet)_m$ is the equalizer of the maps
\begin{equation}\label{equali}
s, t:  \prod_{\phi \in ob\; D_m} 
X_\phi \to \prod_{\alpha \in mor\; D_m}X_\alpha 
\end{equation}
where 
$X_\alpha= X_\phi$ for $\alpha :[\phi]\to [\phi']$, 
and on the component indexed by $\alpha$, $s$ is the projection from 
$X_\phi$,
whereas $t$ is the projection from $X_{\phi'}$ composed with $\alpha^*$.
In particular, we obtain 
$$\cosk_n \stackrel{\sim}{\longrightarrow} \cosk_n\cosk_m$$ for $n\leq m$.
Similarly, for $n\leq m$, the restriction functor $(i^m_n)^*$ from 
$m$-truncated simplicial sets to $n$-truncated simplicial sets
has a left adjoint and a right adjoint $(i^m_n)_*$.
If $n\leq m$, then 
\begin{equation}\label{simp}
\cosk_n \stackrel{\sim}{\longrightarrow} \cosk_m\cosk_n
\end{equation}
(this is wrongly stated in SGA 4 V 7.1.2).
Indeed, we just apply the following to $i^*_n$: 
$$ (i_n)_* = (i_m)_*(i^m_n)_*  =  (i_m)_*i_m^* (i_n)_*.$$

A simplicial map 
$\Delta[1]\times X_\bullet \to Y_\bullet$ can be described as a collection
of maps $h_\tau :X_j\to Y_j$ for every $\tau:[j]\to [1]$, such that
for $\alpha:[i]\to [j]$ one has $\alpha^*h_\tau = h_{\tau\alpha}\alpha^*$. 
A simplicial homotopy between two maps $f,g:X_\bullet\to Y_\bullet$
is a simplicial map $\Delta[1]\times X_\bullet \to Y_\bullet$ 
such that $h_{c_0^j}=f$ and $h_{c_1^j}=g$, where 
$c_\epsilon^j:[j]\to [1]$ is the constant map to $\epsilon\in [1]$.

\begin{lemma}
\label{brianl}\cite[Lemma 3.0.2.4]{sga4}
Let $f: A_\bullet\to B_\bullet$ and $g:B_\bullet \to A_\bullet$ be maps of
simplicial schemes such that 
\begin{enumerate}
\item $f_p:A_p\to B_p$ is inverse to $g_p$ for $p<n$;
\item $g_n$ is a section to $f_n: A_n\to B_n$;
\item $A_\bullet \cong \cosk_n A_\bullet$, and $B_\bullet \cong \cosk_n B_\bullet$.
\end{enumerate}
Then $f$ and $g$ are simplicial homotopy inverse to each other.
\end{lemma}

\proof
It suffices to prove that two simplicial maps $f,g:X_\bullet \to Y_\bullet$
on $n$-truncated schemes which are equal on $n-1$-truncation
induce homotopic maps on the $\cosk_n$. 
For $i\leq n$ and $\tau:[i]\to [1]$, we define $h_\tau=f_i$ if $\tau=c_0$
and $h_\tau=g_i$ otherwise. If easily follows from the first condition that 
for $i,j\leq n$ and $\alpha:[i]\to [j]$ we have 
$\alpha^*h_\tau = h_{\tau\alpha}\alpha^*$. In degrees $p>n$ we define
the map $h_\tau$ for $\tau:[p]\to [1]$ by the following diagram
$$ \begin{CD}
X_p @> h_\tau >>Y_p\\
@| @| \\
\lim_{\phi\in ob D_m'}X_\phi @>h_{\tau\phi}>> \lim_{\phi\in ob D_m'}
Y_\phi
\end{CD}$$
where $X_\phi=X_i$ for $\phi: [i]\to [p]$ and . 
The maps $h_{\tau\phi}:X_i\to Y_i$
are compatible with the inverse system because for $\alpha:[j]\to [i]$
we have $\alpha^*h_{\tau\phi}= h_{\tau\phi\alpha}\alpha^*:X_i\to Y_j$ by
definition of $h$. It is now easy to see that this induces a
map of simplicial objects, i.e. that it is compatible with the 
maps induced by $[q]\to [p]$. On the other hand, taking $h_{c_0^p}$
we recover $(\cosk_n f)_p$, and taking $h_{c_1^p}$ 
we recover $(\cosk_n g)_p$. Indeed, the maps $h_{\tau\phi}$ between 
inverse systems will be the system of maps $f_i$ and $g_i$, respectively.
\proofend

We will apply this in combination with the following lemma:

\begin{lemma}\cite[Lemma 3.0.2.3]{sga4}
Let $f,g:X_\bullet \to Y_\bullet$ simplicial homotopic and 
$\f$ a functor to an abelian category. Then $\f(f)$
and $\f(g)$ are homotopic maps of the associated chain complexes. 
\end{lemma}

We apply the above mostly to the category of schemes over a fixed
scheme $X$. If we want to emphasize the dependence on the base scheme
$X$, we write $\cosk_n(X_\bullet/X)$.

We often identify simplicial or multi-simplicial objects  $A_\bullet$ in
an abelian category with its corresponding chain complex without notice.


By \cite[Lemma 5.8]{svbk}, every proper cdh-cover can be dominated
by a composition of abstract blow-ups.
Similarly, we have

\begin{proposition}\label{refine}
Every proper ldh-cover $Y\to X$ of a noetherian scheme
can be dominated by a composition
$S\to T\to X$, where $S\to T$ is a finite flat map of degree
prime to $l$, and $T\to X$ a composition of abstract blow-ups.
\end{proposition}

\proof
The proof is in the spirit of \cite{ichweilII}, \cite{svbk}, and
\cite[Prop.2.4]{kellyK}. 
We proceed by induction on the dimension of $X$. 
Base-changing with a proper cdh-cover, we can assume that $X$ is 
reduced and integral. Let $\eta$ be a point of $Y$ which maps
to the generic point of $X$ and such that $[k(\eta):k(X)]$
is prime to $l$. Let $\tilde Y$ be the closure of $\eta$ in $Y$;
$\tilde Y$ is generically finite of degree prime to $l$ over $X$. 
By the flatification theorem \cite{rg}, there is a blow-up 
$X'\to X$ with center $Z$  of smaller dimension such that the 
strict transform $Y'\to X'$ of $\tilde Y\to X$ is flat, 
hence finite flat surjective of degree prime to $l$. 
The induction hypothesis applied
to the base change $Y\times_XZ\to Z$ gives a factorization
$S'\to T'\to Z$
whose union with $Y' \to X'\to X$ is  a factorization as required.  
\proofend

\section{Proof of the main theorem}
\subsection{\v Cech covers} 
We assume that $\f $ satisfies the hypothesis of Theorem
\ref{main1} or \ref{main2}.

\begin{proposition}\label{base}
Let $f:X_0\to X$ be a proper cdh-covering. Then the 
augmentation map $\cosk_0(X_0/X)\to X$ induces a quasi-isomorphism on 
$\f (-)$. The same holds for proper ldh-coverings if $\f$ is a sheaf
of $\Z_{(l)}$-modules. 
\end{proposition}

\proof
We give the proof for the ldh-case and complexes, the cdh-case follows by 
erasing parts (b) of the proof.

a) Given an abstract blow-up square \eqref{blbl}, 
if the statement of the proposition
holds for the pull-back to $Z'$, $X'$ and $Z$, then it also holds for $X$.
Indeed, we obtain proper ldh-coverings  $Z_0=Z\times_XX_0\to Z$, 
$Z_0'=Z'\times_XX_0\to Z'$ and $X'_0=X'\times_XX_0\to X'$. 
Since the functor $\cosk_n$ commutes with fiber products,
we obtain on each level an abstract blow-up square upon applying the coskeleton functor. 
Thus we obtain a map of distinguished triangles
$$\begin{CD}
\f (\cosk_0(Z'_0/Z')) @>>> \f (\cosk_0(Z_0/Z)) \oplus \f (\cosk_0(X'_0/X')) @>>>
\f (\cosk_0(X_0/X))\\
@VVV @VVV @VVV \\
\f (Z') @>>> \f (Z) \oplus \f (X') @>>>
\f (X)
\end{CD}$$
and if two maps are quasi-isomorphisms then so is the third.

b) If $p: X'\to X$ is finite flat of degree $d$ prime to $l$, and 
the statement of the theorem holds for the pull-back to $X'$,
then it also holds for $X$. Indeed consider the diagram
$$\begin{CD}
\f (\cosk_0(X'_0/X')) @>p'_*>> 
\f (\cosk_0(X_0/X))\\
@V\cong Vf_*'V @VVf_*V \\
\f (X') @>p_*>> \f (X).
\end{CD}$$ 
By hypothesis, $p_*p^*$ induces 
multiplication by the invertible number $d$,
hence $p_*$ is split surjective on homology.
Since $\cosk_0(X_0/X)\times_XX'  \cong \cosk_0(X'_0/X')$, the pull-back
along $p$ is compatible with all the simplicial structure maps,
hence compatibly split by $\frac{1}{d}p_*$ on each level. This implies that 
$f'_*(p')^*= p^*f_*$, and that
$p'_*(p')^*$ induces multiplication by the invertible number 
$d$ on homology, so that $(p')^*$ is split injective on
homology. 
Finally, since $p_*f_*'= f_*p'_*$ is surjective on homology,
so is $f_*$, and since $f_*' (p')^*= p^*f_*$ is injective on homology, 
so is $f_*$.

c) If $f$ has a section $s:X\to X_0$, then the proposition follows
using the contracting homotopy
$ s\times \id : X_0^{\times p} \to X_0^{\times p+1} $,
where the fiber product is taken over $X$.

d) In general, by Prop.\ref{refine}, 
we can dominate $X_0\to X$ be a sequence $X'\to X$
of abstract blow-ups and finite flat maps of degree prime to $l$.
By a) and b) and induction of the dimension, it suffices to 
prove the theorem after base change to $X'$. But then the map
$X_0\times_XX'\to X'$ has a section induced by the map $X'\to X_0$.
\proofend

\subsection{The general case}
We now give the proof of 
Theorem \ref{main1} and Theorem \ref{main2}
in  the spirit of SGA 4 Vbis \S  3 and \cite{gillet}.
Given the hyperenvelope or $l$-hyperenvelope $X_\bullet\to X$, 
let $X_\bullet^n= \cosk_nX_\bullet$ and consider the sequence
\begin{equation}\label{skeleton}
X_\bullet \stackrel{u_n}{\longrightarrow}   
X_\bullet^n \stackrel{v_n}{\longrightarrow} X_\bullet^{n-1} 
\stackrel{v_{n-1}}{\longrightarrow} \ldots \stackrel{v_2}{\longrightarrow}
X_\bullet^1 \stackrel{v_1}{\longrightarrow} X_\bullet^0  
\stackrel{v_0}{\longrightarrow}X.
\end{equation}
The map $u_n$ is an isomorphism in degrees $\leq n$, so by boundedness 
of $\f(-)$ it
suffices to show that the maps $v_n$ induce quasi-isomorphisms on 
$\f (-)$. The case $v_0$ is Proposition \ref{base}.
By \eqref{simp}, $v_{n}$ satisfies the condition of the following 
proposition for $n\geq 1$:

\begin{proposition}\cite[lemma 3.3.3.2]{sga4} \label{sgadesc}
Let $f: K_\bullet \to L_\bullet $ be a map of simplicial schemes such that 
\begin{enumerate}
\item $K_p\to L_p$ is an isomorphism for $p<n$;
\item $K_n\to L_n$ is a proper cdh-covering (ldh-covering);
\item $K_\bullet \cong \cosk_n K_\bullet$, and 
$L_\bullet \cong \cosk_n L_\bullet$.
\end{enumerate}
Then $\f (K_\bullet)\to \f (L_\bullet)$ is a quasi-isomorphism.
\end{proposition}

\proof  (see \cite[Theorem 7.17]{brian})
Let $[K/L]^p_\bullet$ be the $p$th fiber product of $K_\bullet $ over 
$L_\bullet$. Consider the bisimplicial scheme $Z_{\bullet, \bullet}$
with the $(q+1)$-fold fiber product 
$Z_{p,q}= K_p\times_{L_p} \cdots \times_{L_p} K_p$ in bidegree $(p,q)$
such that the $p$th column $Z_{ p,\bullet}$ is $\cosk_0(K_p/L_p)$, 
and the $q$th row  is $Z_{\bullet,q}$ is $[K/L]^{q+1}_\bullet $. 
In particular, the $p$th column is $K_p$ for $p<n$ by hypothesis.
We have the vertical augmentation 
$\tilde f:Z_{\bullet,\bullet} \to L_\bullet$ induced by 
$f: Z_{\bullet,0} = K_\bullet\to L_\bullet$.
$$\begin{CD} 
@VVV @VVV @VVV \\
K_0\times_{L_0}K_0 @<<< K_1\times_{L_1}K_1 @<<< K_2\times_{L_2}K_2@<<< \cdots \\
@VVV @VVV @VVV \\
K_0 @<<< K_1 @<<< K_2 @<<< \cdots \\
@Vf_0VV @Vf_1VV @Vf_2VV \\
L_0 @<<< L_1 @<<< L_2 @<<< \cdots \\
\end{CD}$$

By Proposition \ref{base} and the following lemma,
we can see column by column 
that $\tilde f$ induces a quasi-isomorphism on $\f (-)$:

\begin{lemma}\cite[Lemma 3.3.3.3]{sga4}
Under the assumptions of Prop. \ref{sgadesc}, all maps
$K_m\to L_m$ are proper cdh-coverings (resp. ldh-coverings).
\end{lemma} 

\proof 
By \eqref{equali}, we have a map of equalizers:
$$\begin{CD}
K_m @>>> \prod_{\phi \in ob\; D_m} K_\phi @>s,t>>  
\prod_{\alpha \in mor\; D_m}K_\alpha \\
@VVV  @VfVV @VfVV \\
L_m@>>> \prod_{\phi \in ob\; D_m} L_\phi @>s,t>>  
\prod_{\alpha \in mor\; D_m}L_\alpha,
\end{CD}$$
where $D$ runs through the injections $\phi:[i]\to [m]$, $i\leq n$ and 
$K_\phi:= K_i$, $L_\phi:= L_i$. 
It suffices to show that the fiber product of the left square
$P= L_m\times_{\prod L_\phi} \prod K_\phi$
has the universal property of the equalizer, because then $K_m\to L_m$ 
is a base-change of a proper cdh-covering (ldh-covering), 
hence is itself a  proper cdh-cover (ldh-covering).

Given a map $u:T\to \prod_{\phi \in D_m} K_\phi $
with $su=tu$, we have to show that there is a unique map $T\to P$ such that
composition with the projection $P\to \prod_{\phi \in D_m} K_\phi$ is $u$.  
By definition of 
$L_m$ and $P$, it suffices to show that $sfu=tfu$, or alternatively that 
$fsu=ftu$, and we can do this factor by factor. 
Given $\alpha:[i]\to [i']$, there are two cases: If $i<n$,
then the two maps agree because $K_i\to L_i$ is an isomorphism by hypothesis.
If $i=n$, then also $i'=n$ and $\alpha$ is the identity, hence $s=t$ trivially.  
\proofend


We apply Lemma \ref{brianl} to $f$ any of the face (i.e. projection)
maps $[K/L]^{p+1}_\bullet\to [K/L]^p_\bullet$, and 
$g$ a degeneracy (i.e. diagonal) map which is a section
of this face map. Note that the hypothesis   
$K_\bullet \cong \cosk_n K_\bullet$, and 
$L_\bullet \cong \cosk_n L_\bullet$
are preserved under fiber products (since $\cosk_n$ is a right adjoint),
hence $[K/L]^p_\bullet\cong \cosk_n [K/L]^p_\bullet$. 
Lemma \ref{brianl} implies that all face maps 
$[K/L]^{p+1}_\bullet\to [K/L]^p_\bullet$ induce quasi-isomorphism
on $\f (-)$ which are equal on homology, 
hence taking alternating sum of projection maps, we see that the 
maps between rows of $Z_{\bullet,\bullet}$ induce alternatingly the 
zero map on homology and quasi-isomorphisms on $\f (-)$. 
This implies that the inclusion
$K_\bullet \to Z_{\bullet,\bullet}$ as the $0$th row induces a 
quasi-isomorphism on $\f (-)$, hence so does the composition
$$f:K_\bullet \to Z_{\bullet,\bullet}\stackrel{\tilde f}{\to}L_\bullet.$$
\proofend

\section{Applications}
\subsection{Existence of $l$-hyperenvelopes}

For $X$ a noetherian scheme, and $l$ a prime number, 
a morphism $h : X' \to  X$  is called an $l$-alteration if
$h$ is proper, surjective, generically finite, sends each maximal point 
to a maximal
point, and the degrees of the residual extensions $k(x')/k(x)$ 
over each maximal point $x$ of $X$ are prime to $l$. 

\begin{theorem}\cite[X Theorem 2.1]{gabber}
(Gabber). Let $k$ be a field, $l$ a prime number different from the
characteristic of $k$, $X$ separated and finite type over $k$.
Then there exists a finite extension $k'$ of $k$ of degree prime to $l$,
and a projective $l$-alteration $h : \tilde X \to X$ above
$\Spec k' \to  \Spec k$, such that $\tilde X$ is smooth
and quasi-projective over $k'$.
\end{theorem}

\begin{corollary}\label{canfindcover}
In the situation of the theorem, there is a proper ldh-cover
over $\tilde X\to X$ with $\tilde X$ regular, and 
quasi-projective, separated and of finite type over $k$.
\end{corollary}

\proof
Covering $X$ by its reduced irreducible components we can assume
that $X$ is integral, and proceed by noetherian induction.
Let $\tilde X$ be as in the theorem, and take $X'$ the closure 
of a point of $\tilde X$ mapping to the generic point of $X$
such that the degree of the residue extension if finite prime to $l$.
Let $Z$ be the closed subscheme where $X'\to X$ is not flat, then
by induction hypothesis we can find a proper ldh-cover $Z'\to Z$
with $Z'$ regular, 
and $Z'\coprod X'\to X$ is the required ldh-covering.
\proofend

To prove Theorem \ref{hyperexist}, we apply the method
of \cite[6.2.5]{deligne} to the Corollary.

\subsection{Motivic theories}
Let $k$ be a perfect field. For a presheaf with transfers $\P$ on $\Sm/k$,
recall that $\underbar C_*(\P)$ 
is the complex of presheaves with transfers given by 
$\underbar C_i(\P)(U)=\P(U\times \Delta^i)$ and boundary maps given by alternating sum of pull-backs along embeddings of faces. 
The complex of abelian groups 
$\underbar C_*(\P)(k)$ is denoted by $C_*(\P)$.

Recall that the cdh-topology is the coarsest Grothendieck topology
generated by Nisnevich covers and proper cdh-covers, and the ldh-topology
is the coarsest Grothendieck topology generated by cdh-covers and
finite flat maps of degree prime to $l$.

\begin{theorem}\label{fvtheorem} 
1) Let $\P$ be a presheaf with transfers such that
$\P_{cdh}=0$. Then under resolution of singularities,
the complex of Nisnevich sheaves with transfers 
$\underbar C_*(\P)_{Nis}$ is acyclic. 

2) If $\P$ is a presheaf of $\Z[\frac{1}{p}]$-modules with transfers such
that $\P_{ldh}=0$, then $\underbar C_*(\P)_{Nis}$ is acyclic. 
\end{theorem}

\proof
This is \cite[Theorem 5.5(2)]{fv}, and its extension by Kelly
\cite[Thm. 5.3.1]{kellythesis}. 
In loc.cit. $\P$ is supposed to be a
presheaf with transfers on all schemes, but replacing $\P$ by its left Kan extension to all schemes does not change $\P$ on smooth schemes. 
Moreover, loc.cit. assumes that $\P_{cdh}=0$,
but in fact for a sheaf of  $\Z[\frac{1}{p}]$-modules,  $\P_{ldh} = 0$ implies that $\P_{cdh} = 0$.
\proofend


Recall that a morphism $f:X\to Y$ is called equidimensional of relative
dimension $r$ if  it is of finite type, if every irreducible component of
$X$ dominates an irreducible component of $Y$, and if 
$\dim_x(p^{-1}p(x))=r$ for every point $x\in X$. A morphism is
equidimensional if and only if it can be locally factored as 
$X\to \mathbb A^r_Y\to Y$, with the first map quasi-finite and
dominant on each irreducible component.
A useful criterion is that a flat morphism of finite type is
equidimensional of dimension $r$ if all 
irreducible components of all generic fibers have dimension $r$.  
For any scheme $X$ over $k$, let $z_{equi}(X,r)$ be the presheaf with
transfers 
on $\Sm/k$ which associates to $U$ the free abelian group on those closed
integral subschemes of $U\times X$ which are equidimensional
of relative dimension $r$ over $U$, and 
$c_{equi}(X,0)$ the subpresheaf with transfers 
of $z_{equi}(X,0)$ generated by 
those subschemes which are finite over $U$.

\begin{proposition}\label{svlemma}
For any scheme $W$, the functors $C_*(z_{equi}(-\times W,r))$ and 
$C_*(c_{equi}(-\times W,0))$ satisfy the hypothesis of Theorems \ref{main1}
and \ref{main2}.
\end{proposition}

\proof
By the existence of transfers, a finite 
flat map $f:X\to Y$ induces a pull-back map such that the composition with push-forward
is multiplication by the degree, and which is compatible with base change.
To show the exactness of the triangle resulting from an abstract
blow-up, it suffices by Theorem \ref{fvtheorem} to show that 
the functors  $\P=z_{equi}(-\times W,r)$ and 
$\P=c_{equi}(-\times W,0)$ send abstract
blow-up squares to short exact sequences of ldh-sheaves (or cdh-sheaves
under resolution of singularities). Replacing $X,X',Z,Z'$ by 
its product with $W$, we can drop $W$ from the notation. Only the 
surjectivity of $\P(X')\oplus \P(Z)\to \P(X)$ is difficult. 

We repeat the platification argument of Suslin-Voevodsky 
\cite[Theorem 4.7]{svbk} and Friedlander-Voevodsky \cite[Theorem 5.11]{fv}. 
Given a section $S\in \P(X)(U)$, we need to find a
cdh-covering (ldh-covering) $V\to U$ such that $S|_V$ is in the image of  
$\P(X')(V)\oplus \P(Z)(V)$. This is clear if $S\subseteq U\times Z$. 
Otherwise let $T$ be the closure of 
$S\cap U\times (X-Z) \subseteq U\times (X-Z)\cong U\times (X'-X'\times_XZ)$
in $U\times X'$. 

Then $T$ may not be equidimensional, but 
by the flatification theorem \cite{rg}, we can find a blow-up
$U'\to U$ such that the proper transform $T'$ of $T$ in $U'\times X'$
is flat over $U'$. 
By Corollary \ref{canfindcover}, we can find a cdh-cover 
(respectively ldh-cover) $V\to U'$  with $V$ smooth, and $T^*$ be
the pull-back of $T'$ to $V$. Then $T^*\to V$ is flat,
hence equidimensional.
\proofend

We define motivic Borel-Moore homology to be
$$ H_i^{BM}(X,A(n)) = 
\begin{cases}
H_{i-2n}C_*(z_{equi}(X,n)\otimes A)& n\geq 0;\\
H_{i-2n}C_*(z_{equi}(X\times \mathbb A^{-n},0)\otimes A) &n<0.
\end{cases}$$
Under resolution of singularities, this agrees with the 
definition of Friedlander-Voevodsky
\cite[\S 9]{fv} by \cite[Thm.5.5(1)]{fv}.
Recall that motivic homology is defined by 
$$ H_i(X,A(n)) = 
\begin{cases}
H^{2n-i}_{\{0\}}
(\mathbb A^n,\underline C_*(c_{equi}(X,0)\otimes A))& n\geq 0;\\
H_{i-2n-1}C_*(\frac{c_{equi}(X\times (\mathbb A^{-n}-\{0\}),0)}
{c_{equi}(X\times \{1\},0)} \otimes A) &n<0.
\end{cases}$$


\medskip

{\it  Proof of theorem} \ref{desccor}: 
For Borel-Moore homology, the theorem follow from Prop. \ref{svlemma}. 
For homology in negative degrees, the theorem follows because the inclusion
$c_{equi}(X\times \{1\},0)\to c_{equi}(X\times (\mathbb A^{-n}-\{0\}),0)$
is canonically split by the structure map, so that the
conclusion of Proposition \ref{svlemma} applies to this situation.

Finally, motivic homology in positive degrees is the homology 
of the cone of the split injection 
$$C_*(c_{equi}(X,0))\otimes A\to R\Gamma(\mathbb A^n,
\underline C_*(c_{equi}(X,0)\otimes A)).$$
Consider the Godement resolution for the Nisnevich topology
of the affine space. We can bound it, because the cohomological
dimension of $\mathbb A^n$ is $-n$.
Furthermore, by \cite[Example 6.20]{mvw}, the terms of
this resolution are still presheaves with transfers. 
Since the Godement resolution is functorial and 
sends short exact sequences of sheaves to short exact sequences of
sheaves, it inherits the hypothesis of Theorems
\ref{main1} and \ref{main2} from  $c_{equi}(X,0)$.
\proofend
 
{\it  Proof of theorem} \ref{desccorchow}:
This follows easily from the localization property of higher 
Chow groups \cite{bloch}\cite{levine}.

\subsection{Other applications}
In \cite{ichsuslin} we show that Parshin's conjecture implies that 
for a smooth variety over a finite field, higher Suslin homology 
vanishes rationally, i.e. 
 $H_p^S(Y,\Q)=0$ for smooth $Y$ and $p>0$.
This implies that for every $l$-hyperenvelope consisting of 
smooth schemes there is an isomorphism
$$ H_p^S(X,\Q) \cong H_p(H_0^S(X_\bullet,\Q)).$$

In \cite{ichrojtman}, we use Corollary \ref{desccor} to show that
for a normal connected variety $X$ over an algebraically closed field,
the Albanese map
$$ alb_X:  H_0(X,\Z)^0\to \Alb_X(k)$$
is an isomorphism on torsion groups away from the characteristic,
and at the characteristic under resolution of singularities.

\end{document}